# IMMERSED PROJECTIVE PLANES, ARF INVARIANTS AND EVEN 4–MANIFOLDS

CHRISTIAN BOHR

ABSTRACT. In this paper, we exploit a subtle indeterminacy in the definition of the spherical Kervaire–Milnor invariant which was discovered by R. Stong to construct non–spin 4–manifolds with even intersection form and prescribed signature.

As an important ingredient of the proof of their $\pi_1$–negligible embedding theorem, M. Freedman and F. Quinn introduced in [4] a $\mathbb{Z}_2$–valued invariant of an immersed sphere representing a spherically characteristic homology class, which is called the spherical Kervaire–Milnor invariant. Later on, R. Stong pointed out that there is a subtle indeterminacy in the definition of the Kervaire–Milnor invariant which is related to the existence of non–trivial immersed projective planes [11]. The main objective of this note is to demonstrate how this indeterminacy can be exploited to construct 4–manifolds with prescribed fundamental group and even intersection form.

Recall that the intersection form of a 4–manifold is called even if the self intersection number of every embedded surface is even, or, equivalently, if there is a characteristic homology class which is torsion. Clearly a 4–dimensional spin manifold has even intersection form. However, the example of the Enriques surface shows that there are 4–manifolds with even intersection form which do not carry a spin structure. In this paper, we shall study the geography of these manifolds, i.e. the question for the pairs of integers which can be realized as Euler characteristic and signature of a non–spin 4–manifold with even intersection form. Note that these two numbers together with the fundamental group completely determine the intersection form of such a manifold.

Clearly the Euler characteristic of a 4–manifold with even intersection form is even and its signature is a multiple of eight. If such a manifold is not spin, its first integral homology group must contain non–trivial 2–torsion. In [12], P. Teichner showed that conversely, every finitely presentable group whose first homology contains 2–torsion is the fundamental group of a 4–manifold with even intersection form which is not spin. Moreover the signature of this manifold can be chosen to be any given multiple of eight. Unfortunately the methods developed in [12] do not seem suitable to study the geography of even 4–manifolds, essentially because the Euler characteristic is not a bordism invariant.

A more geometrical approach to constructing even 4–manifolds is to perform surgery along an embedded 2–sphere whose homology class differs from a characteristic class only by addition of a torsion class. In some cases, such a sphere can be found by exploiting the aforementioned indeterminacy in the definition of the Kervaire–Milnor invariant. As an example of the results which can be obtained in this way, we have the following theorem, which will be proved in section 3.

**Theorem 1.** *Suppose that $\Pi$ is one of the groups $\mathbb{Z}_n$, $\mathbb{Z} \oplus \mathbb{Z}_n$ or $\mathbb{Z}^3 \oplus \mathbb{Z}_n$, where $n$ is even. Assume further that $x$ and $y$ are non–negative integers. If $x$ is even, $y$ is a multiple of eight and $\frac{11}{8} y \leq x - 2$, then there exists a non–spin even 4–manifold with fundamental group $\Pi$, signature $y$ and Euler characteristic $x$.*

The main open conjecture related to the geography of spin 4–manifolds is the 11/8–conjecture, which predicts the inequality
$$\frac{11}{8} |\sigma(X)| \leq b_2(X)$$
between the signature and the second Betti number of a 4–dimensional spin manifold $X$. This conjecture is widely believed to be true, and slightly weaker inequalities were proved by M. Furuta in [5] and by M. Furuta and Y. Kametani in [6].

As demonstrated in [1], the 11/8–conjecture would imply similar inequalities between the signature and the Betti–numbers of a 4–manifold with even intersection form which is not spin. Therefore one should not hope for a comprehensive result on the geography of even 4–manifolds as long as this fundamental conjecture is not settled. However, it seems nevertheless interesting to analyze the consequences it would have. It turns out that assuming this conjecture, one can distill most of the information on the geography of 4–manifolds with even intersection form and abelian fundamental group into a single group invariant.

**Definition 1.** Assume that $\Pi$ is a finitely generated group whose first homology $H_1(\Pi; \mathbb{Z})$ contains non–trivial 2–torsion. Define
$$\tilde{r}(\Pi) = \min\{\chi(X) - \frac{11}{8}|\sigma(X)| \,\big|\, X \text{ even}, w_2(X) \neq 0, \pi_1(X) = \Pi\} \in \mathbb{Z} \cup -\infty,$$
where $\chi(X)$ denotes the Euler characteristic of a 4–manifold $X$.

Unfortunately proving the finiteness of the invariant $\tilde{r}(\Pi)$ is essentially as hard as proving the 11/8–conjecture. In fact, if the conjecture turned out to be true, one could adapt the arguments in [1] to prove that the invariant $\tilde{r}(\Pi)$ is always finite. Conversely the finiteness of $\tilde{r}(\Pi)$ for a group $\Pi$ implies the 11/8–conjecture for simply connected 4–manifolds, as one could take connected sum with sufficiently many copies of a simply connected 4–manifold contradicting the conjecture to make the number $\chi(X) - \frac{11}{8}|\sigma(X)|$ smaller than any given integer. If, however, the invariant *is* finite, it gives an almost complete answer to the geography question in the case of abelian fundamental groups. In fact, we have the following result, which shall be proved in section 3.

**Theorem 2.** *Assume that $\Pi$ is a finitely generated abelian group with non–trivial 2–torsion for which the number $\tilde{r}(\Pi)$ is finite. Then all but finitely many pairs $(x,y)$ of integers for which $x$ is even, $y \equiv 0 \mod 8$ and*
$$\frac{11}{8}|y| \leq x - \tilde{r}(\Pi)$$
*can be realized as $(\chi(X), \sigma(X))$ for a non–spin 4–manifold $X$ with even intersection form and fundamental group $\Pi$.*

It is clear that – assuming again the finiteness of $\tilde{r}(\Pi)$ – these conditions on $x$ and $y$ are also necessary.

Theorem 1 and Theorem 2 are proved using a result concerning embedded spheres in 4–manifolds, which we will now describe. Recall that a homology class $\xi \in H_2(X; \mathbb{Z})$ of a 4–manifold $X$ is called spherically characteristic if $\xi \cdot S \equiv S \cdot S$ mod 2 for every immersed sphere $S$. A standard approach to constructing an embedded sphere is to perform surgery along an embedded surface whose fundamental group is mapped trivially to the fundamental group of the ambient 4–manifold. Following [4], we will call such an embedding a $\pi_1$–null embedding.

Given a $\pi_1$–null embedding $F \to X$ which represents a spherically characteristic homology class $\xi$, there is a $\mathbb{Z}_2$–valued obstruction to reducing the genus of $F$. This invariant is called the Arf invariant $\mathrm{Arf}(F)$. We remark that the usual definition of the Arf invariant of a characteristic surface as described for instance in [3] has to be modified slightly to work for surfaces which are only spherically characteristic, see [2]. If the Arf invariant of $F$ vanishes, one can find an embedded sphere in $X \# k(S^2 \times S^2)$ for some non–negative integer $k$ representing again the class $\xi$.

In the characteristic case, it is known [3] that the Arf invariant is determined by the homology class. This is no longer true if the homology class of $F$ is spherically characteristic, but not characteristic. There are good reasons to believe that in this case, the class $[F]$ can always be stably represented by a surface of Arf invariant zero.

**Conjecture.** *Suppose that $F \to X$ is a $\pi_1$–null embedding representing a spherically characteristic homology class $\xi$ which is not characteristic. Then there is a number $k$ and a $\pi_1$–null embedding $F' \to X \# k(S^2 \times S^2)$ which represents $\xi$ and has Arf invariant zero.*

Using bordism theory, this was proved in [2] with the additional assumption that there exists a spherical class $\omega \in H_2(X; \mathbb{Z})$ such that $\omega \cdot \xi = 1$. In this note, we shall demonstrate that, for certain fundamental groups, one can also find a purely geometric proof based on Stong's observation which does not require the existence of the class $\omega$. Apart from giving further evidence to the above conjecture, this shows that there is a relation between the bordism theory used in [2], which in turn is closely related to earlier work of P. Teichner [12], and the indeterminacy discovered by Stong.

Our results are particularly easy to state if the fundamental group of the manifold at hand is abelian and if the homology class is almost characteristic, in the following sense.

**Definition 2.** *Suppose that $X$ is a 4–manifold. A homology class $\xi \in H_2(X; \mathbb{Z})$ is called almost characteristic if it can be written as $\xi = c + a$ where $c$ is characteristic and $a$ is a torsion class.*

As the pullback of a torsion class to the universal covering is zero, an almost characteristic class is spherically characteristic. However there are spherically characteristic classes which are not almost characteristic. Moreover it is not difficult to see that a 4–manifold has even intersection form if and only if the trivial homology class is almost characteristic. We now have the following result, a more general version of which will be stated and proved in section 2.

**Theorem 3.** *Suppose that $X$ is a 4–manifold and that $F \to X$ is a $\pi_1$–null embedding representing an almost characteristic homology class $\xi$. If $\pi_1(X)$ is abelian and the class $\xi$ is not characteristic, then there is another $\pi_1$–null embedding $F' \to X$ representing $\xi$ such that $\mathrm{Arf}(F') \neq \mathrm{Arf}(F)$. In particular $\xi$ can be stably represented by an embedded sphere.*

Throughout this paper, all manifolds will be assumed to be smooth, oriented, connected and closed, unless stated otherwise. Similarly all maps will be smooth and immersions of surfaces in 4–manifolds will be assumed to be generic in the sense that they only have finitely many singular points which are transverse self intersections.

1. Labeled immersions and Kervaire–Milnor invariants

In this short section, we set up some notation and define a simple modification of the Kervaire–Milnor invariant which avoids the problems described in [11] by fixing a distinguished preimage of every double point.

**Definition 3.** A **labeled immersion** $(f, O)$ is an immersion $f\colon S^2 \to X$ together with a subset $O$ of the set of singular points of $f$ which contains exactly one preimage of every double point.

Suppose that we are given an immersion $f\colon S^2 \to X$ with vanishing Wall self intersection number $\mu(f)$ representing a spherically characteristic homology class. Assume further that the signs of all self intersection points of $f$ add up to zero, which can always be achieved by performing cusp homotopies.

Usually the self intersection number $\mu(f)$ is only defined modulo the subgroup generated by all elements $g - g^{-1}$ of $\mathbb{Z}[\pi_1(X)]$, where $g \in \pi_1(X)$. The reason is that the loop associated to a double point depends on the choice of a first and second sheet. However, given a set $O$ as in the above definition, we can associate a loop to a double point $p$, namely the image of a path on $S^2$ running from the preimage of $p$ in $O$ to the other preimage. Adding up the homotopy classes of all these loops, counted with signs, we obtain a representative of $\mu(f)$ in the group ring $\mathbb{Z}[\pi_1(X)]$. The fact that $\mu(f)$ vanishes implies that we can always choose the set $O$ such that this representative is zero. A labeled immersion for which this is the case will be called *admissible*.

Given an admissible labeled immersion $(f, O)$, we can arrange the self intersection points in pairs $(p_i^+, p_i^-)$ such that all the points $p_i^+$ have sign $+1$, the points $p_i^-$ have sign $-1$ and the group elements given by $p_i^+$ and $p_i^-$ coincide for every $i$. Let $x_i^\pm, y_i^\pm$ denote the preimages of $p_i^\pm$, ordered such that $x_i^\pm \in O$. For every $i$, we can choose a path $\alpha_i$ on $S^2$ from $x_i^-$ to $x_i^+$ and another path $\beta_i$ from $y_i^+$ to $y_i^-$. The images of these two arcs define an oriented loop in $X$. This loop is nullhomotopic and is therefore the boundary of an embedded disk $W_i \subset X$. Following the terminology used in [10], we will call the image of $\alpha_i$ the *positive arc* $\partial_+ W_i$ and the image of $\beta_i$ the *negative arc* $\partial_- W_i$.

The restriction of the normal bundle of $W_i$ to $\partial W_i$ has a non–vanishing section, which is tangential to $f$ along $\partial_+ W_i$ and normal to $f$ along $\partial_- W_i$. As performing a boundary twist changes the obstruction to extending this section by $\pm 1$, we can also assume that this section can be extended over $W_i$ (i.e. that $W_i$ is "correctly framed"). For every Whitney disk $W_i$, we then have an intersection number $W_i \cdot f$ counting the number of intersection points between $f$ and the interior of $W_i$. Now we can define the *Kervaire–Milnor invariant* $km(f, O) \in \mathbb{Z}_2$ of $(f, O)$ by

$$km(f, O) = \sum_i W_i \cdot f \mod 2.$$

The arguments in [4] and [11] now show that this invariant is well defined and only depends on the map $f$ and the set $O$, although it may actually depend on the choice of this set.

Given a $\pi_1$–null embedded surface $F \subset X$ realizing a spherically characteristic homology class, one can define an Arf invariant $\mathrm{Arf}(F) \in \mathbb{Z}_2$ by adapting the definition which is usually used in the characteristic case [2]. The relation between this invariant and the Kervaire–Milnor invariant is provided by the following lemma.

**Lemma 1.** *Suppose that $(f, O)$ is an admissible labeled immersion representing a spherically characteristic homology class. Then there exists a $\pi_1$–null embedded surface $F \to X$ realizing the homology class $[f]$ whose Arf invariant equals $km(f, O)$.*

*Proof.* Assume that we have chosen Whitney disks $W_i$ as above. Remove the images of small disks in $S^2$ around $x_i^\pm$ and join the resulting boundary circles by a 1–handle centered around the negative arc of $W_i$. Performing this for all $i$ results in a $\pi_1$–null embedded surface $F$ realizing the class $[f]$. A symplectic basis for $H_1(F; \mathbb{Z}_2)$ is given by the boundaries of the $W_i$ and the boundaries of the cocores $V_i$ of the attached 1–handles. The cocore of every 1–handle is a correctly framed disk and intersects $f$ transversely in exactly one point. We can use the collection $\{W_i, V_i\}$ to compute the Arf invariant of $F$. As the $W_i$ are correctly framed, we obtain that

$$\mathrm{Arf}(F) = \sum_i q(V_i) q(W_i) = \sum_i f \cdot W_i = km(f, O)$$

as claimed. □

## 2. Detectable classes and projective planes

The normal 1–type of a 4–manifold whose universal covering is spin is determined by a finitely presentable group $\Pi$ together with a distinguished element of $H^2(\Pi; \mathbb{Z}_2)$ [9]. This pair is called the $w_2$–type of the manifold. Similarly a spherically characteristic homology class has a $w_2$–type. We will be interested in a certain property of $w_2$–types, which, roughly speaking, encodes whether the difference between a characteristic class and a class which is only spherically characteristic can be detected by immersed projective planes.

**Definition 4.** A $w_2$–*type* is a pair $(\Pi, w)$, where $\Pi$ is a finitely presentable group and $w \in H^2(\Pi; \mathbb{Z}_2)$. A $w_2$–type $(\Pi, w)$ is called *detectable* if there exists a homomorphism $\varphi \colon \mathbb{Z}_2 \to \Pi$ such that $\varphi^* w \neq 0$.

Although it is not difficult to find examples of $w_2$-types which are not detectable (such an example even exists for groups as simple as $\mathbb{Z}_2 \times \mathbb{Z}_4$), we shall now see that there is an interesting class of detectable $w_2$-types, namely those which correspond to 4–manifolds with even intersection form and abelian fundamental group.

**Lemma 2.** *If $\Pi$ is finitely generated abelian group and*

$$w \in \mathrm{Ext}(\Pi; \mathbb{Z}_2) \subset H^2(\Pi; \mathbb{Z}_2),$$

*a non–zero cohomology class, then $(\Pi, w)$ is detectable.*

*Proof.* By the universal coefficient theorem, we have an exact sequence

$$0 \longrightarrow H_2(\Pi; \mathbb{Z}) \longrightarrow H_2(\Pi; \mathbb{Z}_2) \xrightarrow{\rho} \mathrm{Tor}(\Pi; \mathbb{Z}_2) \longrightarrow 0$$

As, by assumption, $w \in Ext(\Pi; \mathbb{Z}_2)$, the evaluation map $x \mapsto \langle w, x \rangle$ vanishes on the image of $H_2(\Pi; \mathbb{Z})$ and therefore defines a map $\mathrm{Tor}(\Pi; \mathbb{Z}_2) \to \mathbb{Z}_2$. Since $w \neq 0$, there is an element $x \in \text{Tor}(\Pi; \mathbb{Z}_2)$ which is not in the kernel of this map. It is not hard to see that we can find a homomorphism $\varphi \colon \mathbb{Z}_2 \to \Pi$ such that $\text{Tor}(\varphi; \mathbb{Z}_2)$ maps the non–trivial element of $\text{Tor}(\mathbb{Z}_2; \mathbb{Z}_2) = \mathbb{Z}_2$ to $x$. This implies that the image $\varphi_*(t)$ of the non–trivial element $t \in H_2(\mathbb{Z}_2; \mathbb{Z}_2) = \text{Tor}(\mathbb{Z}_2; \mathbb{Z}_2)$ is some class in $H_2(\Pi; \mathbb{Z}_2)$ which is sent to $x$ by $\rho$. Hence

$$\langle \varphi^* w, t \rangle = \langle w, \varphi_* t \rangle = 1,$$

which shows that $\varphi^* w \neq 0$. $\square$

We remark that one can of course construct examples of detectable $w_2$–types $(\Pi, w)$ for which $\Pi$ is not abelian, for instance by taking products of abelian groups and groups with suitable cohomology.

Now suppose that we are given a spherically characteristic homology class $\xi$ of a 4–manifold $X$. Choose a classifying map

$$u \colon X \longrightarrow K(\pi_1(X), 1)$$

for the universal covering $\pi \colon \tilde{X} \to X$. By assumption, the class $w_2(X) - PD(\xi)$ is mapped to zero by $\pi^*$, consequently it is of the form $u^* w$ for a uniquely determined $w \in H^2(\pi_1(X); \mathbb{Z})$.

**Definition 5.** The pair $(\pi_1(X), w)$ is called the $w_2$–type of the class $\xi$.

Note that the $w_2$–type of a spherically characteristic homology class is only well defined up to automorphisms of $\pi_1(X)$. If the universal covering of $X$ is spin, for instance because $X$ has even intersection form, the trivial homology class is spherically characteristic. Its $w_2$–type is usually called the $w_2$–type of the manifold $X$.

Now we are ready to state and prove the main result of this section, which is considerably more general than Theorem 3 in the introduction.

**Theorem 4.** *Suppose that $F \to X$ is a $\pi_1$–null embedding realizing a spherically characteristic class $\xi$. If the $w_2$–type of $\xi$ is detectable, then there exists a $\pi_1$–null embedding $F' \to X$ representing $\xi$ such that $\text{Arf}(F) \neq \text{Arf}(F')$. In particular, $\xi$ can be stably represented by an embedded sphere.*

*Proof.* Let $\Pi = \pi_1(X)$. As $F$ is a $\pi_1$–null embedding, we can choose an immersed sphere $f \colon S^2 \to X$ representing the class $\xi$ which has vanishing Wall self intersection number. Choose a classifying map $u \colon X \to K(\Pi, 1)$ for the universal covering. By assumption, the $w_2$–type $(\Pi, w)$ of $\xi$ is detectable, i.e. there exists a homomorphism $\varphi \colon \mathbb{Z}_2 \to \Pi$ such that $\varphi^* w \neq 0$. Let $a = \varphi(1)$. As usual we will also assume that the signs of all the self intersection points of $f$ add up to zero.

Now we are going to exploit the indeterminacy discovered by Stong. Pick a subset $O$ of the set of singular points of $f$ which contains exactly one preimage of every double point, such that the labeled immersion $(f, O)$ is admissible. Choose Whitney disks pairing up the self intersection points according to $O$, as in the definition of $km(f, O)$. After performing a finger move guided by $a$, we can also assume that there is a Whitney disk $W_i$ with group element $a$, pairing up two self intersection points $p_i^+$ and $p_i^-$. Now let

$$O' = (O \setminus \{x_i^+\}) \cup \{y_i^+\},$$

i.e. we interchange the roles of $x_i^+$ and $y_i^+$. As $a^2 = 1$, the labeled immersion $(f, O')$ is again admissible. Choose a path $\alpha'$ on $S^2$ running from $x_i^-$ to $y_i^+$ and a path $\beta'$ from $x_i^+$ to $y_i^-$. Then the images of $\alpha'$ and $\beta'$ form the boundary of another correctly framed embedded Whitney disk $W_i'$. The union of the arcs $\alpha, \alpha', \beta$ and $\beta'$ is a loop $c$ on $S^2$ which is the boundary of a disk $D$. We can assume that the restriction of $f$ to the interior of $D$ is an embedding. The disks $W$ and $W'$ and the disk $f(D)$ together form an immersed $\mathbb{R}P^2$ $A$ in $X$. Note that the image of $\pi_1(A)$ in $\pi_1(X)$ is generated by $a$. As the disks $W_i$ and $W_i'$ are correctly framed, we have (see also [10])

$$km(f, O) - km(f, O') = f \cdot W_i - f \cdot W_i' = w_2(A) + f \cdot A \in \mathbb{Z}_2.$$

Due to Lemma 1, the proof is finished once we can prove that this number is not zero. Let $\iota \colon \mathbb{R}P^2 \to X$ denote the immersion given by $A$ and consider the composition $u \circ \iota$. By the choice of $a$, the induced map

$$(u \circ \iota)_* \colon \pi_1(\mathbb{R}P^2) = \mathbb{Z}_2 \longrightarrow \Pi$$

is the map $\varphi$. This implies that $\iota^* u^* w \neq 0$. As $u^* w = PD(\xi) - w_2(X)$, we have

$$\langle w_2(X), [A] \rangle - f \cdot A = \langle u^* w, [A] \rangle = \langle \iota^* u^* w, [\mathbb{R}P^2] \rangle.$$

But we have just seen that the element $\iota^* u^* w \in H^2(\mathbb{R}P^2; \mathbb{Z}_2)$ is the orientation. Therefore its evaluation on $[\mathbb{R}P^2]$ is not zero and the theorem is proved. □

*Proof of Theorem 3.* Suppose that $\xi = c + a$ is almost characteristic and let $(\Pi, w)$ denote its $w_2$–type. It is not difficult to show that $w \in \mathrm{Ext}(H_1(\Pi); \mathbb{Z}_2)$. As $\xi$ is not characteristic, $w \neq 0$. Moreover $\Pi = \pi_1(X)$ is abelian by assumption, and therefore the $w_2$–type $(\Pi, w)$ is detectable by Lemma 2. Now the result follows from Theorem 4. □

## 3. Constructing even 4–manifolds

In this section, we use our results on embedded spheres to construct 4–manifolds with given $w_2$–type and given signature. It is not difficult to see that every $w_2$–type can be realized by a 4–manifold, such a manifold can for instance be obtained by surgery along circles in $k(S^1 \times S^3)$ for some natural number $k$. The correct framings can be determined using the language of $B$–structures as in [9]. However, every manifold constructed in this way will have signature zero. In [12], P. Teichner used stable homotopy theory and calculations involving spectral sequences to prove several results on the signatures of 4–manifolds with a given $w_2$–type. In the case of detectable $w_2$–types, we are now in a position to give geometric proofs of some of these facts. Our main technical result is the following.

**Theorem 5.** *Suppose that $X$ is a non–spin 4–manifold with detectable $w_2$-type $(\Pi, w)$. Then there exists a 4–manifold $Y$ having the same $w_2$–type as $X$ such that $\sigma(Y) = \sigma(X) - 8$ and $\pi_1(Y) = \pi_1(X)$.*

*Moreover $Y$ can be chosen such that $Y \# \mathbb{C}P^2 = X \# \mathbb{C}P^2 \# 8\overline{\mathbb{C}P^2} \# 2(S^2 \times S^2)$, in particular $b_2(Y) = b_2(X) + 12$.*

*Proof.* By assumption, there exists a map $\varphi \colon \mathbb{Z}_2 \to \Pi$ such that $\varphi^* w \neq 0$. Let $a = \varphi(1)$.

Now choose an immersion $S^2 \to \mathbb{C}P^2$ with only one self intersection point realizing the class $3\gamma$, where $\gamma$ denotes the class of a generic complex line. Such an immersion can for instance be obtained as a holomorphic curve of degree 3 with one ordinary double point. Perform a cusp homotopy to introduce an additional self intersection point of the opposite sign and a finger move inside $X \# \mathbb{C}P^2$ guided by an arc representing $a$ to obtain an immersion $f \colon S^2 \to X \# \mathbb{C}P^2$ representing $(0, 3\gamma)$.

Then $f$ has four self intersection points, two with group element $a$ and two with trivial group element. Moreover the homology class represented by $f$ is spherically characteristic with $w_2$–type $(\Pi, w)$, in particular its $w_2$–type is detectable.

As in the proof of Theorem 4, we can now choose a labeling $O$ of $f$ such that $km(f, O) = 0$. Using the same construction as in the proof of Lemma 1, we obtain a $\pi_1$–null embedding $F \to X \# \mathbb{C}P^2$ representing the class $(0, 3\gamma)$ such that $F$ has genus 2 and vanishing Arf invariant. Now the arguments given by A. Yasuhara in [13] show that we can replace this surface by an embedded sphere $S \subset X \# \mathbb{C}P^2 \# 2(S^2 \times S^2)$ whose homology class is $[S] = (0, 3\gamma, (0, 2a), (0, 2b))$ for integers $a$ and $b$. Note that this class is spherically characteristic with $w_2$–type $(\Pi, w)$ and has self intersection number 9. Blowing up at eight points on $S$ yields an embedded sphere
$$S' \subset X \# \mathbb{C}P^2 \# 8\overline{\mathbb{C}P^2} \# 2(S^2 \times S^2)$$
which still has $w_2$–type $(\Pi, w)$ such that $S' \cdot S' = 1$. A neighborhood of $S'$ has boundary $S^3$. If we replace that neighborhood by a 4–ball, we obtain a 4–manifold $Y$ which has $w_2$–type $(\Pi, w)$ such that
$$Y \# \mathbb{C}P^2 = X \# \mathbb{C}P^2 \# 8\overline{\mathbb{C}P^2} \# 2(S^2 \times S^2).$$
In particular $Y$ has signature $\sigma(Y) = \sigma(X) - 8$ as desired. □

**Example 1.** Suppose that $\Pi$ is abelian and that $w \in \text{Ext}(\Pi; \mathbb{Z}_2)$ is not zero. Then a 4–manifold with $w_2$–type $(\Pi, w)$ has even intersection form, but is not spin. As already mentioned above, one can easily realize $(\Pi, w)$ by a 4–manifold with signature zero. By Lemma 2, $(\Pi, w)$ is detectable. Therefore Theorem 5 shows that every finitely generated abelian group with non–trivial 2–torsion is the fundamental group of a 4–manifold with even intersection form and signature $-8$, a fact which was proved in [12] using completely different techniques.

**Example 2.** Given an even number $n$, there exists a rational homology 4–sphere with fundamental group $\mathbb{Z}_n$ which is not spin, but has even intersection form. Such a manifold can be obtained by surgery along circles in $S^1 \times S^3$. Applying Theorem 5 shows that there is a non–spin even 4–manifold $X$ such that $\pi_1(X) = \mathbb{Z}_n$, $\sigma(X) = -8$ and $b_2(X) = 12$. Hence the intersection form of $X$ is $2H \oplus E_8$, where $H$ denotes as usual the standard hyperbolic form of rank 2. Guided by the example of the Enriques surface, it is natural to conjecture that $X$ splits off a copy of $S^2 \times S^2$.

Although this example seems to suggest that the construction used in the proof of Theorem 5 can be improved, it is worth noting that such an improvement without further restrictions on the fundamental group would contradict the 11/8–conjecture (this can easily be proved using a covering trick as in [1]). Hence we should only hope to improve Theorem 5 at the cost of restricting to a special class of fundamental groups, for instance to cyclic or – more general – finite groups.

Theorem 5 also has some consequences regarding the geography of 4–manifolds with even intersection form, and, more generally, the geography of 4–manifolds with a given $w_2$–type. To simplify the discussion, we restrict ourselves to oriented 4–manifolds with non–negative signature, keeping in mind that reversing the orientation changes the sign of the signature. Before we can state our results, we first have to introduce the following invariant of $w_2$–types.

**Definition 6.** Suppose that $\Pi$ is a finitely presentable group and that we are given a class $w \in H^2(\Pi; \mathbb{Z}_2)$. Let

$$q_0(\Pi, w) = \min\{\chi(X) \,|\, X \text{ is a 4–manifold with } w_2\text{–type } (\Pi, w), \sigma(X) = 0\}.$$

We remark that this invariant is finite and bounded from below by the invariant $q(\Pi)$ defined by J. Hausmann and S. Weinberger in [7], by the number $p(\Pi)$ introduced by D. Kotschick in [8] and by the invariant $r(\Pi)$ defined in [2].

**Example 3.** Suppose that $n$ is a natural number. Performing surgery along a circle in $S^1 \times S^3$, one can construct a rational homology 4-sphere with fundamental group $\mathbb{Z}_n$. If $n$ is even, this manifold can be chosen to be spin or non–spin, depending on the framing. As every 4–manifold with fundamental group $\mathbb{Z}_n$ must have Euler characteristic at least two, this implies that $q_0(\Pi, w) = 2$ if $\Pi$ is a finite cyclic group.

**Theorem 6.** *Assume that we are given a detectable $w_2$–type $(\Pi, w)$. Suppose that $(x, y)$ is a pair of integers where $y \geq 0$ is a multiple of eight and $x$ is even. If*

$$\frac{11}{8}y \leq x - q_0(\Pi, w),$$

*then there exists a 4–manifold $X$ with $w_2$–type $(\Pi, w)$, signature $\sigma(X) = y$ and Euler characteristic $\chi(X) = x$.*

*Proof.* By the definition of $q_0(\Pi, w)$, there is a 4–manifold $X_0$ with $w_2$–type $(\Pi, w)$ such that $\chi(X_0) = q_0(\Pi, w)$ and $\sigma(X_0) = 0$. First let us suppose that $y$ is actually a multiple of 16. Consider the manifold $X_1$ which is obtained from $X_0$ by adding $\frac{y}{16}$ copies of the K3–surface with the non–complex orientation. Then again the $w_2$–type of $X_1$ is $(\Pi, w)$, $\text{sign}(X_1) = y$ and $\chi(X_1) = q_0(\Pi, w) + \frac{22}{16}y$. By assumption, we have $\chi(X_1) \leq x$, and of course the difference $d = x - \chi(X_1)$ is even. Let $X = \frac{d}{2}(S^2 \times S^2) \# X_1$.

Now assume that $y = 16k + 8$ for some non–negative number $k$. As before we can obtain a 4–manifold $X_1$ with $w_2$–type $(\Pi, w)$, Euler characteristic $q_0(\Pi, w) + 22k$ and signature $-16k$ by taking the connected sum of $X_0$ with $k$ copies of a Kummer surface. By Theorem 5, there exists a 4–manifold $X_2$ with $w_2$–type $(\Pi, w)$ such that

$$\chi(X_2) = \chi(X_1) + 12 = q_0(\Pi, w) + 22k + 12 = q_0(\Pi, w) + \frac{11}{8}y + 1$$

and $\sigma(X_2) = -16k - 8$. By assumption, $\chi(X_2) \leq x + 1$. But $x$ is even and $\chi(X_2)$ is even, hence equality cannot occur. Therefore the difference $d = x - \chi(X_2)$ is again a non–negative even number and we can obtain the desired manifold $X$ by adding $\frac{d}{2}$ copies of $S^2 \times S^2$ to $X_2$ and reversing the orientation. □

*Proof of Theorem 1.* As $\Pi$ contains $\mathbb{Z}_2$, there is a non–zero cohomology class $w \in \text{Ext}(\Pi; \mathbb{Z}_2)$. By Lemma 2, $(\Pi, w)$ is detectable. Every 4–manifold with $w_2$–type $(\Pi, w)$ has even intersection form but is not spin. As in Example 3, one can easily prove that $q_0(\Pi, w) \leq 2$ in all the three cases, using surgery along circles in $S^1 \times S^3$, $T^2 \times S^2$ and $T^4$. Hence the result follows from Theorem 6. □

*Proof of Theorem 2.* Choose a 4–manifold $X_0$ with $\frac{11}{8}|\sigma(X_0)| = \chi(X) - \tilde{r}(\Pi)$ which has even intersection form and fundamental group $\Pi$, but is not spin. We can assume that the signature $\sigma(X_0)$ is not negative. The $w_2$–type of $X_0$ is $(\Pi, w)$ for some non–zero element $w \in \text{Ext}(\Pi; \mathbb{Z}_2)$, in particular it is detectable. Obviously $q_0(\Pi, w) \geq \tilde{r}(\Pi)$.

The region in the $(x, y)$–plane which is cut out by the lines $y = 0$, $y < \sigma(X_0)$, $\frac{11}{8}y = x - \tilde{r}(\Pi)$ and $\frac{11}{8}y = x - q_0(\Pi, w)$ contains only finitely many integral points. Hence our result follows from Theorem 5 once we can show that every point $(x, y)$ with $y \geq \sigma(X_0)$, $y \equiv 0 \mod 8$, $x \equiv 0 \mod 8$ and $\frac{11}{8}y \leq x - \tilde{r}(\Pi)$ can be realized.

Given such a lattice point $(x, y)$, let $d = y - \sigma(X_0)$. By assumption, $d$ is non–negative and a multiple of 8. We only work out the argument in the case that $d \equiv 8 \mod 16$, the case that $d$ is a multiple of 16 is actually easier. Write $d = 16k + 8$. Let $K$ denote the K3–surface with the non–complex orientation and consider the manifold $X_1 = X_0 \# kK$. The $w_2$–type of $X_1$ is again $(\Pi, w)$. Furthermore $\sigma(X_1) = y - 8$ and

$$\chi(X_1) = \chi(X_0) + 22k = \frac{11}{8}\sigma(X_0) + \tilde{r}(\Pi) + \frac{22}{16}(d-8) = \tilde{r}(\Pi) + \frac{11}{8}y - 11.$$

By Theorem 5, we can find a 4–manifold $X_2$ with even intersection form which is not spin, has fundamental group $\Pi$, signature $\sigma(X_1) + 8 = y$ and Euler characteristic

$$\chi(X_2) = \chi(X_1) + 12 = \tilde{r}(\Pi) + \frac{11}{8}y + 1.$$

By assumption, $\frac{11}{8}y + \tilde{r}(\Pi) \leq x$. As $\frac{11}{8}d$ is odd, equality cannot occur, hence $\chi(X_2) \leq x$. Therefore we can obtain the desired manifold by adding $\frac{1}{2}(x - \chi(X_2))$ copies of $S^2 \times S^2$ to $X_2$. □

Mathematisches Institut, Theresienstr. 39, 80333 München, Germany
*E-mail address*: bohr@mathematik.uni-muenchen.de